\documentclass[a4paper, 12pt]{article}
\usepackage{makeidx}
\usepackage{enumerate, theorem}
\usepackage{amsmath, amsfonts, amssymb}
\usepackage[height=22.5cm, width=15cm]{geometry}
\usepackage[all]{xy}
\usepackage{mathrsfs, yfonts}
\usepackage{graphicx}
\newcommand{\parag}[1]{\paragraph{\sc{#1.}} }

\DeclareMathOperator{\Sym}{Sym}

\newtheorem{thm}{Theorem}[subsection]

\newtheorem{cor}[thm]{Corollary}
\newtheorem{prop}[thm]{Proposition}
\newtheorem{lemma}[thm]{Lemma}

\begin{document}

\title{Note on rational $1-$dimensional compact cycles. second version}
\date{27/09/16}

 \author{Daniel Barlet\footnote{Institut Elie Cartan, G\'eom\`{e}trie,\newline
Universit\'e de Lorraine, CNRS UMR 7502   and  Institut Universitaire de France.}.}

 \maketitle
 
 \parag{Abstract}
 We give a proof of the fact tha the subset of the rational curves  form a closed analytic subset in the space of  the $1-$dimensional cycles of a complex space.\\
 
 \parag{Key words} Rational curves - Cycles' space.\\
 
 \parag{AMS Classification 2010} 32 C 25 - 32 G 10.\\
 
 The aim of this short note is to prove the following result.
 
 \begin{thm}\label{rat. 0}
 Let $M$ be a complex space. Let $R$ be the subset of the space $\mathcal{C}_{1}(M)$ consisting of compact rational $1-$dimensional cycles in $M$. Then $R$ is a closed analytic subset in $\mathcal{C}_{1}(M)$.
 \end{thm}
 
 By a rational $1-$dimensional cycle we mean that each irreducible component of such a cycle is a rational curve (may be singular).\\
 
  Note that this result is classical in the projective context.\\
  
  The proof of the theorem uses the following proposition.
  
  \begin{prop}\label{closed}
  Let $\pi : U \to V$ be a geometrically flat map between reduced complex spaces with one dimensional  fibres. Assume that for a point $v_{0} \in V$ the fibre $\pi^{-1}(v_{0})$ has an irreducible component $\gamma_{0}$ which has genus $\geq 1$. Then there exists an open neighbourhood $V_{1}$ of $v_{0}$ in $V$ such that for any $v \in V_{1}$ the fibre $\pi^{-1}(v)$ has an irreducible component of genus $\geq 1$.
   \end{prop}
  
  In the previous statement, the genus of an irreducible compact curve is, by definition, the genus of its normalization.\\
  Recall also that the geometric flatness assumption means that there exists a holomorphic map $\varphi :V \to \mathcal{C}_{1}(U)$ such that for each $v \in V$ we have $\pi^{-1}(v)  = \vert \varphi(v)\vert$ and such that for $v$ generic in $V$ the $1-$cycle $\varphi(v)$ is reduced (that is to say that all multiplicities are equal to $1$).

  We begin by some general results in order to show that it is enough to prove the proposition and the theorem in the case when $V$ is normal and when the  fibres of $\pi$ are connected.
  
  \begin{lemma}\label{base modif.}
  Let $U \to V$ be a geometrically flat map between reduced complex spaces with holomorphic  fibre map $\varphi : V \to \mathcal{C}_{n}(U)$. Let $\tau : \tilde{V} \to V$ be a proper modification of $V$ and denote $\tilde{\pi} : \tilde{U} \to \tilde{V}$ the strict transform of $\pi$ by $\tau$. Then $\tilde{\pi}$ is geometrically flat and the fibre at a point $\tilde{v} \in \tilde{V}$ is the $n-$cycle $\varphi(\tau(\tilde{v}))\times \{\tilde{v}\}$ in $\tilde{U}$.
  \end{lemma}
  
  \parag{proof} Let $S \subset V$ be the center of the modification $\tau$. Then for $\tilde{v} \not\in \tau^{-1}(S)$ the fibre of $\tilde{\pi}$ is $\varphi(\tau(\tilde{v}))\times \{\tilde{v}\}$ as a cycle in $\tilde{U} \subset U \times_{V} \tilde{V}$. As $\tilde{U}$ is a closed analytic subset in  $U \times_{V} \tilde{V}$ and as the family of cycles $\tilde{v} \mapsto \varphi(\tau(\tilde{v}))\times \{\tilde{v}\}$ is an analytic family of compact cycles in $U \times_{V} \tilde{V}$ such that, for $\tilde{v}$ generic, they are contained in $\tilde{U}$, this is an analytic  family of cycles in $\tilde{U}$ and this gives the holomorphic fibre map for $\tilde{\pi}$.$\hfill \blacksquare$\\
  
  \parag{Remarks} \begin{enumerate}
  \item If  $n = 1$ and if $R $ is the set of points in $\tilde{V}$ such that the fibre of $\tilde{\pi}$ is rational then $\tau(R)$ is the subset of $V$ where the fibre of $\pi$ is rational.
  \item This lemma allows to assume that $V$ is a normal complex space in the proofs of the proposition \ref{closed} and the theorem \ref{rat. 0}
  \end{enumerate}
  
  \begin{lemma}\label{connected fibres}
   Let $U \to V$ be a geometrically flat map between irreducible complex spaces and assume that $V$ is normal. Let $\pi' : U \to W$ and $g : W \to V$ be a Stein factorization of $\pi$. Let $\nu : \tilde{W} \to W$ be the normalization of $W$ and let $ \tilde{\pi} : \tilde{U} \to \tilde{W}$ be the strict transform of $\pi'$ by $\nu$. Put $\tilde{g} : \tilde{W} \to V$ be the composition $\tilde{g} := \nu\circ g$. As $\tilde{\pi}$ is equi-dimensional proper and with a normal basis, it is geometrically flat with connected fibres. Also $\tilde{g}$ is proper, finite and surjective with a normal basis so it is geometrically flat with $0-$dimensional fibres.\\
   Let $\theta : \tilde{W} \to \mathcal{C}_{n}(\tilde{U})$ and $f : V \to \Sym^{k}(\tilde{W})$ the corresponding holomorphic fibre maps. Then the fibre map for the geometrically flat map $ \tilde{\pi}\circ \tilde{g}$ is given by
    $$f\circ \Sym^{k}(\theta)\circ Add : V \to \mathcal{C}_{n}(\tilde{U}),$$
    where $Add : \Sym^{k}(\mathcal{C}_{n}(\tilde{U})) \to \mathcal{C}_{n}(\tilde{U})$ is the addition map of $n-$cycles in $\tilde{U}$.
   \end{lemma}
   
   \parag{proof} Remember that, by definition of a Stein reduction, the fibre of $g$  of a point $v \in V$ is the set of connected components of $\pi^{-1}(v)$ and that the fibre of $\pi'$ at a point $w \in g^{-1}(v)$ is the connected component of $\pi^{-1}(v)$ given by $w$. So $\pi'$ and also $\tilde{\pi}$ have connected fibres. As $g$ is proper (and finite) and as $W$ is irreducible the image by $g$ of non normal points in $W$ is a closed analytic subset in $V$ with no interior point. So it is clear that the holomorphic map $ f\circ \Sym^{k}(\theta)\circ Add$ is a fibre map for $\tilde{\pi}\circ \tilde{g}$ at the generic points of $V$. This is enough to conclude.$\hfill \blacksquare$\\

   \begin{cor}\label{conn. fibre}
   In the situation of the previous lemma with $n = 1$, let $\tilde{R}$ be the subset of $\tilde{W}$ of points such the fibre of $\tilde{\pi}$ is rational. Then the subset of point in $V$ such that the fibre of $\pi$ is rational is equal to $R := \{ v \in V \ / \  f(v) \in \Sym^{k}(\tilde{R})\}$. Then, if $\tilde{R}$ is closed (resp. analytic), so is $R$.
   \end{cor}
   
   \parag{proof} This corollary is clear because a  compact curve is rational if and only if each of its  connected component is rational.$\hfill \blacksquare$\\
   
 \parag{remark}  With the previous results, it is enough to prove the proposition  \ref{closed} and the theorem \ref{rat. 0} with the following extra assumptions :  $V$ is normal and all fibres of $\pi$ are connected.\\

    \parag{proof of the proposition} We shall use the following result  (see [B.80]):  Let $C$ be a reduced compact curve in a complex space $M$, and let $L$ be a holomorphic line bundle on $C$. Then there exists an open neighbourhood $M'$ of $C$ in $M$ and a holomorphic line bundle $\mathcal{L}$ on $M'$ inducing $L$ on $C$. Moreover, if $L$ is topologically trivial on $C$ we may choose $\mathcal{L}$ topologically trivial on $M'$.\\
    Note that the last point is not stated in {\it loc. cit} but is a trivial consequence of the  proof given there.\\

  Consider now  $C_{0} = \vert \varphi(v_{0})\vert = \pi^{-1}(v_{0})$ and let $\nu : \tilde{C}_{0} \to C_{0}$ be the normalization of $C_{0}$. Define the coherent sheaf  $\mathcal{F} := \nu_{*}(\mathcal{O}_{\tilde{C}_{0}})$ on $C_{0}$. We have an exact sequence of coherent sheaves
  $$  0 \to \mathcal{O}_{C_{0}} \to \mathcal{F} \to Q \to 0  $$
  where $Q$ has support in a finite set. Then $H^{1}(C_{0}, Q)$ vanishes and we have a surjective map
  $$ H^{1}(C_{0}, \mathcal{O}_{C_{0}}) \to H^{1}(C_{0}, \mathcal{F}) \to 0 . $$
  Using the surjectivity above, we can find a topologically trivial line bundle $L$ on $C_{0}$ which is not holomorphically trivial on each non rational component of $C_{0}$. Thanks to the result quoted above and to the properness of $\pi$ we can find an open neighbourhood $V_{0}$ of $v_{0}$ in $V$  an  a line bundle $\mathcal{L}$ on $U_{0} := \pi^{-1}(V_{0})$ which is topologically trivial on $U_{0}$ and induces $L$ on $C_{0}$. Now let $Z$ be the subset of the space $\mathcal{C}_{1}(\mathcal{L})$ of connected\footnote{see th. IV 7.2.1 for the analyticity of this condition on compact cycles.} compact $1-$cycles in $\mathcal{L}$ such that their direct image on $U_{0}$ is contained in $\varphi(V_{0}) \subset \mathcal{C}_{1}(U_{0})$. This subset $Z$ is a closed analytic subset in $\mathcal{C}_{1}(\mathcal{L})$ because $\varphi(V_{0})$ is a closed analytic subset in $\mathcal{C}_{1}(U_{0})$ as $\varphi : V_{0} \to \mathcal{C}_{1}(U_{0})$ is a proper  holomorphic map and as the direct image by the projection $p : \mathcal{L} \to U_{0}$ is holomorphic.\\
    Remark that for each $v \in V_{0}$ the set $Z$ contains the cycle $\varphi(v)$ of $\mathcal{L}$ which is the zero section of the restriction of  $\mathcal{L}$ to $\vert \varphi(v)\vert$ with suitable multiplicities, in order that its direct image on $U_{0}$ is equal to $\varphi(v)$. This defines a closed holomorphic embedding of $\varphi(V_{0})$ in $Z \subset \mathcal{C}_{1}(\mathcal{L})$.\\
  We shall show now that the direct image map $ f_{*} : Z \to \varphi(V_{0})$ for compact $1-$cycles induced by the projection $\mathcal{C}_{1}(\mathcal{L}) \to \mathcal{C}_{1}(U_{0})$ has positive dimensional fiber at $\varphi(v) \in Z$ when $\varphi(v)$ is rational.\\
 Assume that for some $v \in V_{0}$ the $1-$cycle $C := \varphi(v)$ is rational. Then the restriction of the line bundle $\mathcal{L}$ on $C$ is holomorphically  trivial and any compact $1-$dimensional cycle in $\mathcal{L}_{\vert C}$ can be move (by vertical translation). So any point in $Z$ in the fibre of $f_{*}$ over $C = \varphi(v)$ is not isolated. This proves our assertion.\\
  Now $C_{0}$ has at least one irreducible component, say $\gamma$, which is not rational; so the corresponding  point of $\varphi(v_{0}) \in Z$  is isolated in its fibre for $f_{*}$. Indeed, the zero section  is the only reduced compact $1-$dimensional cycle in $L_{\vert \gamma} $ as $L_{\vert \gamma}$ is topologically trivial but not holomorphically trivial by construction. So any connected compact $1-$cycle near-by $\varphi(v_{0})$ in $Z\cap f_{*}^{-1}(C_{0})$ must have support in  the zero section of $L$  on each non rational component of $C_{0}$. As $\varphi(v_{0})$  is connected  this implies that on a rational component of $\varphi(v_{0})$ which meets an irrational component, the corresponding component of a near-by cycle to the cycle $\varphi(v_{0})$ in $f_{*}^{-1}(C_{0})$  has to vanish at some point (the intersection with some non rational component). Then the corresponding component of such a cycle is the zero section over this rational component (we have only constant sections on rational components). As we assume $C_{0}$ connected and as the cycles in  $Z$ are connected, we conclude that $\varphi(v_{0}) \in Z$ is an isolated point in its fibre of $f_{*}$.\\
   Now the subset $T$ of points $t \in Z$ such that  the dimension at $t$ of the fibre of $f_{*}$ is at least equal to $1$ is a closed analytic set in $Z$. The intersection of $T$ with the closed embedding of $\varphi(V_{0})$ is $Z$  defines a closed analytic subset in $\varphi(V_{0})$ and then also of $V_{0}$ which contains the subset of rational fibers of $\pi$ in $V_{0}$. As we have shown that $v_{0}$ is not in this closed analytic subset, we obtain an open neighourhood $V_{1}$ of $v_{0}$ in $V_{0}$ such that for any $v \in V_{1}$ the cycle $\varphi(v)$ is not rational. $\hfill \blacksquare$\\

  \begin{cor}\label{rat.1} Let $M := U \times \mathbb{P}_{1}$ be a reduced complex space and let $p: M \to \mathbb{P}_{1}$ be the projection. Let $p_{*} : \mathcal{C}_{1}(M) \to \mathcal{C}_{1}(\mathbb{P}_{1}) \simeq \mathbb{N}.[\mathbb{P}_{1}]$ the direct image for compact $1-$cycles. Define $X := p_{*}^{-1}(1.[\mathbb{P}_{1}])$. As $p_{*}$ is holomorphic, this is a closed analytic subset in $\mathcal{C}_{1}(M)$.
  Let $S_{0}$ be the subset of $\mathcal{C}_{1}(M)$ of irreducible cycles which are in $X$. This is a Zariski open subset in $X$ (see [B-M] prop. IV 7.1.2). Then the closure $S_{1}$ of $S_{0}$ in $X$ contains only rational cycles.
   \end{cor}
   
   \parag{proof} Let $ \pi : U \to V := S_{1}$ be the projection of the graph of the tautological family of $1-$cycles parametrized by $S_{1}$. As the generic cycle in this family is irreducible, by definition of $S_{0}$, all fibres of $\pi$ are connected. Also the generic fibres are rational because for $s \in S_{0}$ it is reduced and isomorphic to $\mathbb{P}_{1}$. If there is a non rational cycle in $S_{1}$, then there exists, thanks to the previous proposition, a non empty open set of non rational cycles in this family. But $S_{0}$ is open and dense in $S_{1}$ this gives a contradiction.$\hfill \blacksquare$
   
   \parag{proof of the theorem} Let $V$ be an irreducible component of $\mathcal{C}_{1}(M)$. As we may normalize $V$ thanks to the lemma \ref{base modif.}, we can assume that the generic cycle in $V$ is reduced. Let $W$ be a relatively compact open set in $V$. Then there exists, thanks to the proposition IV 7.1.2 in [B-M 1], an integer $k \geq 1$ such that for any $v \in W$ the cycle $v$ has at most $k$ irreducible components.\\
   Denote $\pi : U \to W$ the projection of the graph of the tautological family of compact curves in $M$ parametrized by $W$.\\
    Define the subset $S_{l} \subset \mathcal{C}_{1}(U \times \mathbb{P}_{1})$ as the image by the addition map of cycles of $(S_{1})^{l}$ in $\mathcal{C}_{1}(U \times \mathbb{P}_{1})$. As the addition map is proper and finite, $S_{l}$ is a closed analytic subset of $\mathcal{C}_{1}(U \times \mathbb{P}_{1})$ for each integer $l \geq 1$.\\
   Let $q : U \times \mathbb{P}_{1} \to U$ the projection. We shall prove that the subset of rational cycles in the family $(v)_{v \in W}$ is exactly given by the subset
   $$ R := \big(\cup_{l = 1}^{k} \ q_{*}(S_{l}) \big) \cap W$$
   which is a closed analytic subset in $W$ because the map $q$ is proper and the subset $S_{l} \subset \mathcal{C}_{1}(U \times \mathbb{P}_{1})$ is a closed analytic subset.\\
   First remark that each cycle in $R$ is rational as the direct image of a rational cycle is rational and we proved that each cycle in $S_{1}$ is rational in corollary \ref{rat.1}.\\
   Conversely, let $v \in W$ such that $v$ is rational. Then there exists an integer $l \in [1, k]$ and $l$ holomorphic generically injective\footnote{the normalization maps.} maps (distinct or not) $f_{1}, \dots, f_{l}$ from $\mathbb{P}_{1}$ to $U$ such that the sum of there images is $v$ and with graphs $G_{1}, \dots, G_{l}$ which are points in $S_{0}$. Then $v = \sum_{i=1}^{l} q_{*}(G_{i})$ and so $v$ is in $R$.\\
   To conclude the proof it is enough to say that a  subset which is closed and analytic on any open relatively compact subset in $V$ is a closed analytic subset in $V$.$\hfill \blacksquare$
   
   \parag{References}
\begin{itemize}
\item{ [B. 78]}  Barlet, D. \ {\it Majoration du volume ...} in  Lecture Notes $n^{0}$ 822 (1980) (Sem. Lelong-Skoda 78-79), Springer.
\item {[B-M 1]} Barlet, D. et Magnusson, J. \  {\it Cycles analytiques complexes I : th\'eor\`{e}mes de pr\'eparation des cycles},  Cours sp\'ecialis\'es $n^{0} 22$ \\ Soci\'et\'e Math\'ematique de France (2014).
\end{itemize}

   \end{document}